\documentclass[11pt,article]{amsart}
\usepackage{amssymb}
\usepackage{amsfonts}
\usepackage{amsbsy}
\usepackage{latexsym}
\theoremstyle{plain}

 \theoremstyle{remark} 

\newcommand{\Z}{{\mathbb Z}}

\newcommand{\inner}[1]{\langle #1 \rangle}

\newtheorem {theo} {\bf Theorem} [section]
\newtheorem {prop} [theo] {\bf Proposition}
\newtheorem {coro} [theo] {\bf Corollary}
\newtheorem {lem} [theo] {\bf Lemma}
\newtheorem {defi} {\bf Definition}[section]
\newtheorem{exam} {\bf Example}[section]
\newtheorem{rem}{\bf Remark}[section]

\newtheorem*{acknowledgements}{Acknowledgements}
\newcommand{\R}{\mathbb R\,}

\newcommand{\ip}[2]{\langle\,#1\, , \, #2 \, \rangle}

\numberwithin{equation}{section}
\begin{document}
\title[Duality Properties] {A Duality principle for groups}
\author{Dorin Dutkay}\address{Department of Mathematics\\
University of Central Florida\\ Orlando, FL 32816}
\email{ddutkay@mail.ucf.edu}
\author{Deguang Han}\address{Department of Mathematics\\
University of Central Florida\\ Orlando, FL 32816}
\email{dhan@pegasus.cc.ucf.edu}
\author{David Larson}\address{Department of Mathematics\\
Texas A\&M University\\ College Station, TX}
\email{larson@math.tamu.edu}
\date{\today}

\begin{abstract}  The duality principle for Gabor frames states that a
Gabor sequence obtained by a time-frequency lattice is a frame for
$L^{2}(\R^{d})$ if and only if the associated adjoint Gabor
sequence is a Riesz sequence. We prove that this duality principle
extends to any dual pairs of projective unitary representations of
countable groups. We examine the existence problem of dual pairs
and establish some connection with classification problems for
$\rm II_{1}$ factors.  While in general such a pair may not exist for
some groups, we show that such a dual pair always exists for every
subrepresentation of the left regular unitary representation when
$G$ is an abelian infinite countable group or an amenable ICC
group. For free groups with finitely many generators, the
existence problem of such a dual pair is equivalent to the well
known problem about the classification of free group von Neumann
algebras.
\end{abstract}
\maketitle

\section{Introduction}
\setcounter{figure}{0} \setcounter{equation}{0}

Motivated by the duality principle for Gabor representations in
time-frequency analysis  we establish a general duality theory for
frame representations of infinite countable groups, and build its
connection with the classification problem \cite{Conn} of $\rm II_{1}$
factors. We start by recalling some notations and  definitions
about frames.

A {\it frame} \cite{DS} for a Hilbert space $H$ is a sequence
 $\{x_{n}\}$ in $H$ with the property that there
exist positive constants $A, B > 0$ such that
\begin{equation} \label{eg1.1}
A\|x\|^{2}\leq \sum_{g\in G}|\ip{x}{x_{n}}|^{2} \leq B \|x\|^{2}
\end{equation}
holds for every $x\in H$. {\it A tight frame} refers to the case
when $A = B$, and a {\it Parseval frame}  refers to the case when
$A = B = 1$. In the case that (\ref{eg1.1}) hold only for all
$x\in \overline{span}\{x_{n}\}$, then we say that $\{x_{n}\}$ is a
{\it frame sequence}, i.e., it is a frame for its closed linear
span. If we only require the right-hand side of the inequality
(\ref{eg1.1}), then $\{x_{n}\}$ is called a {\it Bessel sequence}.

One of the well studied classes of frames is the class of Gabor
(or Weyl-Heisenberg) frames: Let $\mathcal{K}= A\Z^{d}$ and
$\mathcal{L}= B\Z^{d}$ be two full-rank lattices in $\R^d$, and
let $g \in L^2(\R^d)$ and $\Lambda = \mathcal{L}\times
\mathcal{K}$. Then the {\it Gabor (or Weyl-Heisenberg ) family} is
the following collection of functions in $L^2(\R^d)$:
$$
   {\bf G}(g, \Lambda) =  {\bf G}(g, \mathcal{L}, \mathcal{K}):= \{e^{2\pi i<\ell, x>} g(x-\kappa) \Bigm| \ell\in\mathcal{L},
\kappa\in\mathcal{K}\}.
$$
For convenience,  we write $g_{\lambda} = g_{\kappa, \ell} =
e^{2\pi i<\ell, x>} g(x-\kappa)$, where $\lambda = (\kappa,
\ell)$. If $E_{\ell}$ and $T_{\kappa}$ are the modulation and
translation unitary operators defined by $$ E_{\ell}f(x) = e^{2\pi
i<\ell, x>}f(x) $$ and
$$
T_{\kappa}f(x) = f( x - \kappa)
$$
for all $f\in L^{2}(\R ^{d})$. Then we have  $g_{\kappa, \ell} =
E_{\ell}T_{\kappa}g$. The well-known Ron-Shen duality principle
states that a Gabor sequence ${\bf G}(g, \Lambda)$ is a frame
(respectively, Parseval frame) for $L^{2}(\Bbb R^{d})$ if and only
if the adjoint Gabor sequence ${\bf G}(g, \Lambda^{o})$ is a Riesz
sequence (respectively, orthonormal sequence), where $\Lambda^{o}
= (B^t)^{-1}\Z^{d} \times (A^t)^{-1}\Z^{d}$ is the adjoint lattice of
$\Lambda$.

Gabor frames can be viewed as frames obtained by projective
unitary representations of the abelian group $\Z^{d}\times
\Z^{d}$. Let $\Lambda = A\Z^{d} \times B\Z^{d}$ with $A$ and $B$
being $d\times d$ invertible real matrices. The {\it Gabor
representation} $\pi_{\Lambda}$ defined by $(m, n) \rightarrow
E_{Am}T_{Bn}$ is not necessarily a unitary representation of the
group $\Z^{d}\times \Z^{d}$. But it is a projective unitary
representation of $\Z^{d}\times \Z^{d}$. Recall (cf. \cite{Va})
that a {\it projective unitary representation} $\pi$ for a
countable  group $G$ is a mapping $g\rightarrow \pi(g)$ from $G$
into the group $U(H)$ of all the unitary operators on a separable
Hilbert space $H$ such that $\pi(g)\pi(h) = \mu(g, h)\pi(gh)$ for
all $g, h\in G$, where $\mu(g, h)$ is a scalar-valued function on
$G\times G$ taking values in the circle group $\mathbb{T}$. This
function $\mu(g, h)$ is then called a {\it multiplier or
$2$-cocycle} of $\pi$. In this case we also say that $\pi$ is a
$\mu$-projective unitary representation.  It is clear from the
definition that we have

(i) $\mu(g_{1}, g_{2}g_{3}) \mu(g_{2}, g_{3}) = \mu(g_{1}g_{2},
g_{3}) \mu(g_{1}, g_{2})$ for all $g_{1}, g_{2}, g_{3}\in G$,

(ii) $\mu(g, e) = \mu(e, g) = 1$ for all $g\in G$, where $e$
denotes the group unit of $G$.

Any function $\mu: G\times G \rightarrow \mathbb{T}$ satisfying
$(i)$ -- $(ii)$ above will be called a {\it multiplier} for $G$.
It follows from $(i)$ and $(ii)$ that we also have

(iii) $\mu(g, g^{-1}) = \mu(g^{-1}, g)$ holds for all $g\in G$.

Examples of projective unitary representations include unitary
group representations (i.e., $\mu \equiv 1$) and the Gabor
representations in time-frequency analysis.

Similar to the group unitary representation case,   the left and
right regular projective representations with a prescribed
multiplier $\mu$ for $G$  can be defined by
$$
\lambda_{g}\chi_{h} = \mu(g, h)\chi_{gh}, \ \ \  h\in G,
$$
and
$$
\rho_{g}\chi_{h} = \mu(h, g^{-1})\chi_{hg^{-1}}, \ \ \  h\in G,
$$
where $\{\chi_{g}: g\in G\}$ is the standard orthonormal basis for
$\ell^{2}(G)$. Clearly, $\lambda_{g}$ and $r_{g}$ are unitary
operators on $\ell^{2}(G)$. Moreover,  $\lambda$ is a
$\mu$-projective unitary representation of $G$ with multiplier
$\mu(g, h)$ and $\rho$ is a projective unitary representation of
$G$ with multiplier $\overline{\mu(g, h)}$. The representations
$\lambda$ and $\rho$ are called the {\it left regular
$\mu$-projective representation} and the {\it right regular
$\mu$-projective representation}, respectively, of $G$. Let
$\mathcal{L}$ and $\mathcal{R}$ be the von Neumann algebras
generated by $\lambda$ and $\rho$, respectively.  It is known (cf.
\cite{GH1}), similarly to the case for regular group
representations, that both $\mathcal{R}$  and $\mathcal{L}$ are
finite von Neumann algebras, and that $\mathcal{R}$ is the
commutant of $\mathcal{L}$. Moreover, if for each $e \neq u\in G$,
either $\{vuv^{-1}: v\in G\}$ or $\{\mu(vuv^{-1},
v)\overline{\mu(v, u)}: \ v\in G\}$ is an infinite set, then both
$\mathcal{L}$ and $\mathcal{R}$ are factor von Neumann algebras.

{\bf Notations.} In this paper for a subset $M$ of a Hilbert space $H$ and  a
subset $\mathcal A$ of $B(H)$ of all the bounded linear operators
on $H$, we will use $[M]$ to denote the closed linear span  of
$M$,  and ${\mathcal A}'$ to denote the commutant $\{T\in B(H): TA
= AY, \forall A\in\mathcal A\}$ of $\mathcal A$. So we have
$\mathcal{L} = \lambda(G)'' = \rho(G)'$ and $\mathcal{R} =
\rho(G)'' = \lambda(G)'$. We also use
$\mathcal{M}\simeq\mathcal{N}$ to denote two $*$-isomorphic von
Neumann algebras $\mathcal{M}$ and
 $\mathcal{N}$.
\smallskip

Given a projective unitary representation $\pi$ of a countable
group $G$ on a Hilbert space $H$, a vector $\xi \in H$ is called a
{\it complete frame vector (resp. complete tight frame vector,
complete Parseval frame vector)} for $\pi$ if $ \{\pi(g)\xi\}_{
g\in G}$ (here we view this as a sequence indexed by $G$) is a
frame (resp. tight frame, Parseval frame) for the whole Hilbert
space $H$, and is just called a {\it frame vector (resp. tight
frame vector,  Parseval frame vector)} for $\pi$ if $
\{\pi(g)\xi\}_{ g\in G}$  is a {\it frame sequence (resp. tight
frame sequence, Parseval frame sequence)}. A {\it Bessel vector}
for $\pi$ is a vector $\xi\in H$ such that $ \{\pi(g)\xi\}_{ g\in
G}$ is Bessel. We will use $\mathcal{B}_{\pi}$ to denote the set
of all the Bessel vectors of $\pi$.

For $x\in H$, let $\Theta_{x}$ be the analysis operator for $
\{\pi(g)x\}_{ g\in G}$ (see section 2). It is useful to note that
if $\xi$ and $\eta$ are Bessel vectors for $\pi$, then
$\Theta_{\eta}^{*}\Theta_{\xi}$ commutes with $\pi(G)$. Thus, if
$\xi$ is a complete frame vector for $\pi$, then $\eta : =
S_{\xi}^{-1/2}\xi$ is a complete Parseval frame vector for $\pi$,
where $S_{\xi} = \Theta_{\xi}^{*}\Theta_{\xi}$ and is called the
{\it frame operator} for $\xi$ (or {\it Bessel operator} if $\xi$
is a Bessel vector). Hence, a projective unitary representation
has a complete frame vector if and only if it has a complete
Parseval frame vector. In this paper the terminology {\it frame
representation} refers to a projective unitary representation that
admits a complete frame vector.

\begin{prop}\label{prop-sub} \cite{GH1, Rie} \ {\it Let $\pi$ be a  projective unitary representation $\pi$ of a countable
group $G$ on a Hilbert space $H$. Then $\pi$ is frame
representation if and only if $\pi$ is unitarily equivalent to a
subrepresentation of the left regular projective unitary
representation of $G$. Consequently, if  $\pi$ is frame
representation, then both $\pi(G)'$ and $\pi(G)''$ are finite von
Neumann algebras.}
\end{prop}

The duality principle for Gabor frames was independent and
essentially simultaneous discovered by Daubechies, H. Landau, and
Z. Landau \cite{DLL}, Janssen \cite{Jan}, and Ron and Shen
\cite{RS}, and the techniques used in these three articles to
prove the duality principle are completely different. We refer to
\cite{He07} for more details about this principle and its
important applications. For Gabor representations, let
$\Lambda^{o}$ be the adjoint lattice of a lattice $\Lambda$. The
well-known density theorem (c.f. \cite{HW1}) implies that one of
two projective unitary representations $\pi_{\Lambda}$ and
$\pi_{\Lambda^{o}}$ for the group $G = \Z^{d}\times\Z^{d}$ must
be a frame representation and the other admits a Riesz vector. So
we can always assume that $\pi_{\Lambda}$ is a frame
representation of $\Z^{d}\times \Z^{d}$ and hence
$\pi_{\Lambda^{(o)}}$ admits a Riesz vector.  Moreover, we also
have $\pi_{\Lambda}(G)' = \pi_{\Lambda^{(o)}}(G)''$, and both
representations share the same Bessel vectors. Rephrasing the
duality principle in terms of Gabor representations, it states
that $\{\pi_{\Lambda}(m, n)g\}_{m, n\in \Z^{d}}$  is a frame for
$L^{2}(\R^{d})$ if and only if $\{\pi_{\Lambda^{(o)}}(m, n)g\}_{m,
n\in \Z^{d}}$ is a Riesz sequence. Our first main result reveals
that this duality principle  is not accidental and in fact it is a
general principle for any commuting pairs of projective unitary
representations.

\begin{defi} Two projective unitary representations $\pi$ and
$\sigma$ of a countable group $G$ on the same Hilbert space $H$
are called a {\it commuting pair} if $\pi(G)' = \sigma(G)''$.
\end{defi}

\begin{theo}\label{thm-main1} {\it Let $\pi$ be a frame representation and  $(\pi, \sigma)$ be a commuting
pair of projective unitary representations of $G$ on $H$ such that
$\pi$ has a complete frame vector which is also a Bessel vector
for $\sigma$. Then \vspace{2mm}

(i) $\{\pi(g)\xi\}_{g\in G}$ is a frame sequence if and only if
$\{\sigma(g)\xi\}_{g\in G}$ is a frame sequence, \vspace{2mm}

(ii) if, in addition, assuming that $\sigma$ admits a Riesz
sequence, then $\{\pi(g)\xi\}_{g\in G}$ is a frame (respectively,
a Parseval frame) for $H$ if and only if $\{\sigma(g)\xi\}_{g\in
G}$ is a Riesz sequence (respectively, an orthonormal sequence.}
\end{theo}

For a frame representation $\pi$, we will call $(\pi, \sigma)$ a
{\it dual pair} if $(\pi, \sigma)$ is a commuting pair such that
$\pi$ has a complete frame vector which is also a Bessel vector
for $\sigma$,  and $\sigma$ admits a Riesz sequence. We remark
that this duality property is not symmetric for $\pi$ and $\sigma$
since $\pi$ is assumed to be a frame representation and $\sigma$
in general is not. Theorem \ref{thm-main1} naturally leads to the
following existence problem:
\medskip

\noindent{\bf Problem 1.} {\it Let $G$ be a infinite countable
group and $\mu$ be a multiplier for $G$. Does every
$\mu$-projective frame representation  $\pi$ of $G$ admit a dual
pair $(\pi, \sigma)$?}
\medskip

While we maybe able to answer this problem for some special
classes of groups, this is in general open due to its connections
(See Theorem \ref{thm-main2}) with the classification problem of
$\rm II_{1}$ factors which is one of the big problems in von
Neumann algebra theory. It has been a longstanding unsolved
problem to decide whether the factors obtained from the free
groups with $n$ and $m$ generators respectively are isomorphic if
$n$ is not equal to $m$ with both  $n,m >1$. This problem was one
of the inspirations for Voiculescu's theory of free probability.
Recall that the fundamental group $F(\mathcal {M})$ of a type $\rm
II_{1}$ factor $\mathcal{M}$  is an invariant that was considered
by Murray and von Neumann in connection with their notion of
continuous dimension  in \cite{Mvon}, where they proved that that
$F(\mathcal {M}) = \mathbb{R}_{+}^{*}$ when $\mathcal {M}$ is
isomorphic to a hyperfinite type $\rm II_{1}$ factor, and more
generally when it splits off such a factor. For free groups
$\mathcal{F}_{n}$ of $n$-generators, by using Voiculescu's free
probability theory \cite{VN}, Radulescu \cite{Rad-Inv} showed that the
fundamental group $F(\mathcal{M}) = \R_{+}^{*}$ for $\mathcal{M} =
\lambda(\mathcal{F}_{\infty})'$. But the problem of calculating
$F(\mathcal {M})$ for $M=\lambda(\mathcal{F}_{n})'$ with $2\leq
n<\infty$ remains a central open problem in the classification of
$\rm II_{1}$ factors, and it can be rephrased as:
\medskip

\noindent{\bf Problem 2.} {\it Let $\mathcal{F}_{n}$ ($n > 1$) be
the free group of $n$-generators and
$P\in\lambda(\mathcal{F}_{n})'$ is a nontrivial projection.  Is
$\lambda(\mathcal{F}_{n})'$ $*$-isomorphic to
$P\lambda(\mathcal{F}_{n})'P$? }
\bigskip

It is proved in \cite{Rad-Inv} that either all the von Neumann
algebras $P\lambda(\mathcal{F}_{n})'P$ ( $0\neq
P\in\lambda(\mathcal{F}_{n})'$) are  $*$-isomorphic, or no two of them are $*$-isomorphic. Our second main result
established the equivalence of these two problems for free groups.

\begin{theo}\label{thm-ICC} {\it Let $\pi = \lambda_{P}$ be a subrepresentation of the left
regular representation of an ICC (infinite conjugate class) group
$G$ and $P\in\lambda(G)'$ be a projection. Then the following are
equivalent:

(i)  $\lambda(G)'$ and $P\lambda(G)'P$ are isomorphic von Neumann
algebras;

(ii) there exists a group representation $\sigma$ such that $(\pi,
\sigma)$ form a dual pair. }
\end{theo}

The above theorem implies that  the answer to Problem 1 is
negative in general, but  is affirmative  for amenable ICC groups.

\begin{theo}\label{thm-main2} {\it Let $G$ be a countable group and
$\lambda$ be its left regular unitary representation (i.e. $\mu
\equiv 1$). Then we have

(i) If $G$ is either an abelian group or an amenable ICC group,
then for every projection $0\neq P\in \lambda(G)'$, there exists a
unitary representation $\sigma$ of $G$ such that $(\lambda|_{P},
\sigma)$ is a  dual pair.

(iii) There exist ICC groups (e.g., $G =\Z^{2}\rtimes SL(2, \Z)$), such that none of the nontrivial
subrepresentations $\lambda|_{P}$ admits a dual pair.}
\end{theo}

We will give the proof of  Theorem \ref{thm-main1} in section 2
and the proof requires some resent work by the present authors
including the results on parameterizations and dilations of frame
vectors \cite{GH3, Han2_Parseval, HL}, and some result results on
the ``duality properties" for $\pi$-orthogonal and $\pi$-weakly
equivalent vectors \cite{HL_BLM}. The proofs for Theorem
\ref{thm-ICC} and Theorem \ref{thm-main2} will provided in section
3, and additionally we will also discuss some concrete examples
including the subspace  duality principle for Gabor
representations.

\section{The Duality Principle}

We need a series of preparations in order to prove Theorem
\ref{thm-main1}.

For any projective representation $\pi$ of a countable group $G$
on a Hilbert space $H$  and  $x\in H$,  the {\it analysis operator
} $\Theta_{x}$ for $x$ from $\mathcal{D}(\Theta_{x}) (\subseteq
H)$ to $\ell^{2}(G)$ is defined by
$$
\Theta_{x}(y) = \sum_{g\in G}\inner{y, \pi(g)x}\chi_{g},
$$
where $\mathcal{D}(\Theta_{x}) = \{y\in H: \sum_{g\in G}|\inner{y,
\pi(g)x}|^{2} < \infty\}$ is the domain space of $\Theta_{x}$.
Clearly, $\mathcal{B}_{\pi} \subseteq \mathcal{D}(\Theta_{x})$
holds for every $x\in H$. In the case that $\mathcal{B}_{\pi}$ is
dense in $H$, we have that $\Theta_{x}$ is a densely defined and
closable linear operator from $\mathcal{B}_\pi$ to $\ell^2(G)$
(cf. \cite{GH}). Moreover, $x\in \mathcal{B}_{\pi}$ if and only if
$\Theta_{x}$ is a bounded linear operator on $H$, which in turn is
equivalent to the condition that $\mathcal{D}(\Theta_{x}) = H$.

\begin{lem}  \label{lem1.2} \cite{GH} {\it Let $\pi$ be a projective representation
of a countable group $G$ on a Hilbert space $H$ such that
$\mathcal{B}_{\pi}$ is dense in $H$. Then for any  $x\in H$, there
exists $\xi\in\mathcal{B}_{\pi}$ such that}

(i) {\it $\{\pi(g)\xi: g\in G\}$ is a Parseval frame for
$[\pi(G)x]$;}

(ii) {\it $\Theta_{\xi}(H) = [\Theta_{x}(\mathcal{B}_{\pi})]$.}
\end{lem}

\begin{lem} \label{lem1.4} {\it Assume that $\pi$ is a projective
representation of a countable group $G$ on a Hilbert space $H$
such that $\pi$ admits a Riesz sequence and $\mathcal{B}_{\pi}$ is
dense in $H$. If $[\Theta_{\xi}(H)]\neq \ell^{2}(G)$, then there
exists $0\neq x\in H$ such that $[\Theta_{x}(H)]\perp
[\Theta_{\xi}(H)]$.}
\end{lem}
\begin{proof} Assume that $\{\pi(g)\eta\}_{g\in G}$ is a Riesz
sequence. Then we have that $\Theta_{\eta}(H) = \ell^{2}(G)$ and
$\Theta_{\eta}$ is invertible when restricted to $[\pi(G)\eta]$.
Let $P$ be the orthogonal projection from $\ell^{2}(G)$ onto
$[\Theta_{\xi}(H)]$.  Then $P\in \lambda(G)'$ and $P \neq I$. Let
$x = \theta_{\eta}^{-1}P^{\perp}\chi_{e}$. Then $x\neq 0$ and
$[\Theta_{x}(H)]\perp [\Theta_{\xi}(H)]$.
\end{proof}

\begin{lem} \cite{GH, Han2_Parseval} \label{dilation} {\it Assume that $\pi$ is a projective
representation of a countable group $G$ on a Hilbert space $H$
such that $\pi$ admits a complete frame vector $\xi$. If
$\{\pi(g)\eta\}_{g\in G}$ is a frame sequence, then there exists a
vector $h\in H$ such that $\eta$ and $h$ are $\pi$-orthogonal and
$\{\pi(g)(\eta+h)\}_{g\in G}$ is a frame for $H$.}
\end{lem}

Two other concepts are needed.

\begin{defi} Suppose $\pi$ is a projective unitary representation of  a
countable group $G$ on a separable Hilbert space $H$  such that
that the set $\mathcal B_\pi$ of Bessel vectors for $\pi$ is dense in H.
 We will say that two vectors $x$
and $y$ in $H$ are {\it $\pi$-orthogonal} if the ranges of
$\Theta_x$ and $\Theta_y$ are orthogonal, and that they are {\it
$\pi$-weakly equivalent} if the closures of the ranges of
$\Theta_x$ and $\Theta_y$ are the same.
\end{defi}

The following result obtained in \cite{HL_BLM} characterizes the
$\pi$-orthogonality and $\pi$-weakly equivalence in terms of the
commutant of $\pi(G)$.

\begin{lem}\label{Main-HL-08} {\it Let $\pi$ be a projective representation of a
countable group $G$ on a Hilbert space $H$ such that
$\mathcal{B}_{\pi}$ is dense in $H$, and let $x, y\in H$. Then}

(i) {\it $x$ and $y$ are $\pi$-orthogonal if and only if
$[\pi(G)'x] \perp [\pi(G)'y]$ (or equivalently, $x \perp
\pi(G)'y$);}

(ii) {\it $x$ and $y$ are $\pi$-weakly equivalent  if and only if
$[\pi(G)'x] = [\pi(G)'y],$}
\end{lem}

We also need the following parameterization result \cite{GH3,
Han2_Parseval, HL}.

\begin{lem}\label{param} {\it Let $\pi$ be a projective representation of a
countable group $G$ on a Hilbert space $H$ and
$\{\pi(g)\xi\}_{g\in G}$ is a Parseval frame for $H$. Then

(i) $\{\pi(g)\eta\}_{g\in G}$ is a Parseval frame for $H$ if and
only if there is a unitary operator $U\in \pi(G)''$ such that
$\eta = U\xi$;

(ii) $\{\pi(g)\eta\}_{g\in G}$ is a frame for $H$ if and only if
there is an invertible operator $U\in \pi(G)''$ such that $\eta =
U\xi$;

(iii) $\{\pi(g)\eta\}_{g\in G}$ is a Bessel sequence if and only
if there is an operator $U\in \pi(G)''$ such that $\eta = U\xi$,
i.e., $\mathcal{B}_{\pi} = \pi(G)''\xi$.}
\end{lem}

As a consequence of Lemma \ref{param} we have

\begin{coro} \label{coro1} {\it Let $\pi$ be a frame representation of a
countable group $G$ on a Hilbert space $H$. Then

(i) $\mathcal{B}_{\pi}$ is  dense in $H$,

(ii) $\pi$ has a complete frame vector which is also a Bessel
vector for $\sigma$ if and only if $\mathcal{B}_{\pi} \subseteq
\mathcal{B}_{\sigma}$.}

\end{coro}
\begin{proof} (i) follows immediately from Lemma \ref{param}(iii).

For (ii), assume that $\{\pi(g)\xi\}_{g\in G}$ is a frame for $H$
and $\{\sigma(g)\xi\}_{g\in G}$ is also Bessel. Then for every
$\eta\in\mathcal{B}_{\pi}$, we have by Lemma \ref{param} (iii)
there is $A\in \pi(G)''$ such that $\eta = A\xi$. Thus
$\{\sigma(g)\eta\}_{g\in G} = A\{\sigma(g)\xi\}_{g\in G}$ is
Bessel, and so $\eta\in\mathcal{B}_{\sigma}$. Therefore we get
$\mathcal{B}_{\pi} \subseteq \mathcal{B}_{\sigma}$. The other
direction is trivial.
\end{proof}

Now we are ready to prove Theorem \ref{thm-main1}. We divide the
proof into two propositions.

\begin{prop}\label{prop1} {\it Let $\pi$ be a frame representation and $(\pi, \sigma)$ be a commuting
pair of projective unitary representations of $G$ on $H$ such that
$\pi$ has a complete frame vector which is also a Bessel vector
for $\sigma$. Then $\{\pi(g)\xi\}_{g\in G}$ is a frame sequence
(respectively, a Parseval frame sequence) if and only if
$\{\sigma(g)\xi\}_{g\in G}$ is a frame sequence (respectively, a
Parseval frame sequence).} \vspace{2mm}
\end{prop}

\begin{proof} ``$\Rightarrow$ :" Assume that $\{\pi(g)\xi\}_{g\in G}$ is a frame
sequence.  Since  $\pi$ is a frame representation, by the dilation
result (Lemma \ref{dilation}), there exists $h\in H$ such that
$(\xi, h)$ are $\pi$-orthogonal and  $\{\pi(g)(\xi +h)\}_{g\in G}$
is a frame for $H$. If we prove that $\{\sigma(g)(\xi
+h)\}_{g\in G}$ is a frame sequence, then $\{\sigma(g)\xi\}_{g\in
G}$ is a frame sequence. In fact, using the
$\pi$-orthogonality of $\xi$ and $h$ and  Lemma \ref{Main-HL-08},
we get that $[\pi(G)'\xi]\perp [\pi(G)'h]$, and hence
$[\sigma(G)\xi]\perp [\sigma(G)h]$ since $\sigma(G)'' = \pi(G)'$.
Therefore, projecting $\{\sigma(g)(\xi +h)\}_{g\in G}$ onto
$[\sigma(G)\xi]$ we get that $\{\sigma(g)\xi\}_{g\in G}$ is a
frame sequence as claimed. Thus, without losing the generality, we
can assume that $\{\pi(g)\xi\}_{g\in G}$ is a frame for $H$.

By Corollary \ref{coro1}, we have $\xi \in \mathcal{B}_{\pi}
\subseteq \mathcal{B}_{\sigma}$.  From Lemma \ref{lem1.2} we can
choose $\eta\in [\sigma(G)\xi]=: M$ such that $\xi$ and $\eta$ are
$\sigma$-weakly equivalent and $\{\sigma(g)\eta\}_{g\in G}$ is a
Parseval frame for $[\sigma(G)\xi]$. By the parameterizition
theorem (Lemma \ref{param}) there exists an operator $A\in
\sigma(G)''|_{M}$ such that $\xi = A\eta$. Assume that $C$ is the
lower frame bound for $\{\pi(g)\xi\}_{g\in G}$. Then for every
$x\in M$ we have

\begin{eqnarray*}
||x||^{2} & \leq  & {1\over C} \sum_{g\in G}|\inner{x, \
\pi(g)\xi}|^{2} = {1\over C} \sum_{g\in G}|\inner{x, \
\pi(g)A\eta}|^{2}\\
& = & {1\over C} \sum_{g\in G}|\inner{A^{*}x, \ \pi(g)\eta}|^{2} =
{1\over C}||A^{*}x||^{2}.
\end{eqnarray*}
Thus $A^{*}$ is bounded from below and therefore it is invertible
since  $\sigma(G)''|_{M}$ is a finite von Neumann algebra
(Proposition \ref{prop-sub}). This implies that $A$ is invertible
(on $M$) and so $\{\sigma(g)\xi\}_{g\in G} (=
\{A\pi(g)\eta\}_{g\in G}$ is a frame for $M$. \vspace{2mm}

``$\Leftarrow$:" Assume that $\{\sigma(g)\xi\}_{g\in G}$ is a
frame sequence. Applying Lemma \ref{lem1.2} again there exists
$\eta\in [\pi(G)\xi]$ such that $\eta$ and $\xi$ are $\pi$-weakly
equivalent, and $\{\pi(g)\eta\}_{g\in G}$ is a Parseval frame for
$[\pi(G)\xi]$. Using the converse statement proved above, we get that
$\{\sigma(g)\eta\}_{g\in G}$ is a frame sequence for $M :=
[\sigma(G)\eta]$. Since $\xi$ are $\pi$-weakly equivalent, we have
by Lemma \ref{Main-HL-08} that $[\pi(G)'\xi] = [\pi(G)'\eta]$ and
so $M = [\sigma(G)\eta] = [\sigma(G)\xi]$. Thus
$\{\sigma(g)\eta\}_{g\in G}$ is a frame for $[\sigma(G)\xi]$. By
the parameterization theorem (Lemma \ref{param}), there exists an
invertible operator operator $A\in \sigma(G)''|_{M}$ such that
$\xi = A\eta$. Extending $A$ to an invertible operator $B$ in
$\sigma(G)''$, we have $A\eta = B\eta$, and so
$$
\pi(g)\xi = \pi(g)A\eta = \pi(g)B\eta = B\pi(g)\eta.
$$
Thus $\{\pi(g)\xi\}_{g\in G}$ is a frame sequence since
$\{\pi(g)\eta\}_{g\in G}$ is a frame sequence and $B$ is bounded
invertible.

For the Parseval frame sequence case, all the operators $A$ and
$B$ involved in the parameterization are unitary operators and the
rest of the argument is identical to the frame sequence case.
\end{proof}

\begin{prop}\label{prop2} {\it  Let $\pi$ be a frame representation of $G$ on $H$. Assume that $(\pi, \sigma)$ is a commuting
pair of projective unitary representations of $G$ on $H$ such that
such that $\pi$ has a complete frame vector which is also a Bessel
vector for $\sigma$. If $\sigma$ admits a Riesz sequence, then

(i) $\{\pi(g)\xi\}_{g\in G}$ is a frame for $H$ if and only if
$\{\sigma(g)\xi\}_{g\in G}$ is a Riesz sequence.

(ii) $\{\pi(g)\xi\}_{g\in G}$ is a Parseval frame for $H$ if and
only if $\{\sigma(g)\xi\}_{g\in G}$ is an orthonormal sequence.}
\vspace{2mm}
\end{prop}
\begin{proof} (i) ``$\Rightarrow:$" Assume that $\{\pi(g)\xi\}_{g\in G}$ is a
frame for $H$. Then from Proposition \ref{prop1} we have that
$\{\sigma(g)\xi\}_{g\in G}$ is a frame sequence.

Thus, in order to show that $\{\sigma(g)\xi\}_{g\in G}$ is a Riesz
sequence, it suffices to show that $[\Theta_{\sigma, \xi}(H)] =
\ell^{2}(G)$, where $\Theta_{\sigma, \xi}$ is the analysis
operator of $\{\sigma(g)\xi\}_{g\in G}$. We prove this by
contradiction.

Assume that $[\Theta_{\sigma, \xi}(H)] \neq \ell^{2}(G)$. Then, by
Lemma \ref{lem1.4},  there is a vector $0\neq x\in H$ such that
$\Theta_{\sigma, x}(H) \perp \Theta_{\sigma, \xi}(H)$. Since
$\mathcal{B}_{\sigma}$ is dense in $H$ (recall that $\mathcal B_\pi$ is dense in $H$ since $\pi$ is a frame representation), we get by Lemma
\ref{Main-HL-08} that $[\sigma(G)'x] \perp [\sigma(G)'\xi]$ and so
$[\pi(G)x] \perp [\pi(G)\xi]$ since $\sigma(G)' = \pi(G)''$. On
the other hand, since $\{\pi(g)\xi\}_{g\in G}$ is a frame for $H$,
we have $[\pi(G)\xi] = H$ and so we have  $x = 0$, a contradiction.

``$\Leftarrow:$" Assume that $\{\sigma(g)\xi\}_{g\in G}$ is a
Riesz sequence. Then, again by Proposition \ref{prop1} we
$\{\pi(g)\xi\}_{g\in G}$ is a frame sequence. So we only need to
show that $[\pi(G)\xi] = H$.

Let $\eta \perp [\pi(G)\xi]$. So we have $[\pi(G)\eta] \perp
[\pi(G)\xi]$. By Lemma \ref{Main-HL-08}, we have that
$\Theta_{\sigma, \eta}(H) \perp \Theta_{\sigma, \xi}(H)$. But
$\Theta_{\sigma, \xi}(H) = \ell^{2}(G)$ since
 $\{\sigma(g)\xi\}_{g\in G}$ is a Riesz sequence. This implies that $\eta =
 0$, and so $[\pi(G)\xi] = H$, as claimed.

 (ii) Replace ``frame" by ``Parseval frame", and ``Riesz" by
 ``orthonormal", the rest is exactly the same as in (i).
\end{proof}

\section{The Existence Problem}

We will divide the proof of Theorem \ref{thm-main2} into two
cases: The abelian group case and the ICC group case. We deal the
abelian group first, and start with an simple example when $G =
Z$.


\begin{exam}\label{exam1} \ Consider the unitary representation of $\Z$ defined
by $\pi(n) = M_{e^{2\pi i nt}}$ on the Hilbert space $L^{2}[0,
1/2]$. Then $\sigma (n) =  M_{e^{2\pi i 2nt}}$ is another unitary
representation of $\Z$ on $L^{2}[0, 1/2]$.  Note that
$\{\sigma(n)1_{[0, 1/2]}\}_{n\in \Z}$ is an orthogonal basis for
$L^{2}[0, 1/2]$. We have that $\sigma(Z)''$ is maximal abelian and
hence  $\sigma(Z)''= \mathcal{M}_{\infty} = \pi(Z)'$. Moreover a
function $f\in L^{2}[0, 1/2]$ is a Bessel vector for $\pi$
(respectively, $\sigma$) if and only if $f\in L^{\infty}[0, 1/2]$.
So $\pi$ and $\sigma$ share the same Bessel vectors. Therefore
$(\pi, \sigma)$ is a commuting pair with the property that
$\mathcal{B}_{\pi} = \mathcal{B}_{\sigma}$, and $\sigma$ admits a
Riesz sequence.
\end{exam}

It turns out the this example is generic for abelian countable
discrete group.
 \vspace{2mm}

\begin{prop}\label{abeian-thm-main} {\it Let $\pi$ be a unitary frame representation of
an abelian infinite countable discrete group $G$ on $H$. Then there
exists a group representation $\sigma$ such that  $(\pi, \sigma)$
is a dual pair.}
\end{prop}

\begin{proof} Let $\hat{G}$ be the dual group of $G$. Then
$\hat{G}$ is a compact space. Let $\mu$ be the unique Haar measure
of $\hat{G}$. Any frame representation $\pi$ of $G$ is unitarily
equivalent to a representation of the form: $g\rightarrow
e_{g}|_{E}$, where $E$ is a measurable subset of $\hat{G}$ with
positive measure, and $e_{g}$ is  defined by $e_{g}(\chi)= <g,
\chi>$ for all $\chi\in \hat{G}$. So without losing the
generality, we can assume that $\pi(g) = e_{g}|_{E}$.

Let $\nu(F) := {1\over \mu(E)}\mu(F)$ for any measurable subset
$F$ of $E$. Then both $\mu$ and $\nu$ are Borel probability measures
without any atoms. Hence (see \cite{Di}) there is a measure preserving bijection
$\psi$ from $E$ onto $\hat{G}$. Define a unitary representation
$\sigma$ of $G$ on $L^{2}(E)$ by
$$
\sigma(g)f(\chi) = e_g(\psi(\chi))f(\chi), \ \ f\in L^{2}(E).
$$
Then by the  same arguments as in Example \ref{exam1} we have that
$\{\sigma(g)1_{E}\}_{g\in G}$ is an orthogonal basis for
$L^{2}(E)$, and $(\pi, \sigma)$ satisfies all the requirements of
this theorem.
\end{proof}

%
%
\bigskip

\noindent{\bf Proof of Theorem \ref{thm-ICC}}

``(i) $\Rightarrow$ (ii):"  Let $\Phi: \lambda(G)' \rightarrow
P\lambda(G)'P$ be an isomorphism between the two von Neumann
algebras. Note that $tr(A) = <A\chi_{e}, \chi_{e}>$ is a
normalized normal trace for $\lambda(G)'$. Define $\tau$ on
$\lambda(G)'$ by
$$
\tau(A) = {1\over tr(P)} tr(\Phi(A)), \ \ \forall A\in
\lambda(G)'.
$$
Then $\tau$ is also an normalized normal trace for $\lambda(G)'$.
Thus $\tau(\cdot) = tr(\cdot)$ since $\lambda(G)'$ is a factor von
Neumann algebra. In particular we have that
$$
{1\over tr(P)} tr(\Phi(\rho_{g})) = \tau(\rho_{g}) = tr(\rho_{g})
= \delta_{g, e}.
$$
Therefore, if we define $\sigma(g) = \Phi(\rho_{g})$, then
$\sigma$ is a unitary representation of $G$ such that $\sigma(G)''=P\lambda(G)'P=(\lambda(G)P)'=
= \pi(G)'$ and $\sigma$ admits an orthogonal sequence
$\{\sigma(g)\xi\}_{g\in G}$, where $\xi = P\chi_{e}$. Moreover,
for any $A\in \pi(G)''$ we have that $\sigma(g)A\xi =
A\sigma(g)\xi$ and so $A\xi$ is a Bessel vector for $\sigma$. By
Lemma \ref{param} (iii), we know that $\mathcal{B}_{\pi} =
\pi(G)''\xi$. Thus we get $\mathcal{B}_{\pi}\subseteq
\mathcal{B}_{\sigma}$ and therefore $(\pi, \sigma)$ is a dual
pair.

``(ii) $\Rightarrow$ (i):" Assume that $(\pi, \sigma)$ is a dual
pair. Let $\{\sigma(g)\psi\}_{g\in G}$ be a Riesz sequence,
$\sigma_{1}(g) := \sigma(g)|_{M}$ and $\sigma_{2}(g) :=
\sigma(g)|_{M^{\perp}}$, where $M = [\sigma(g)\psi]$. Then
$\sigma$ is unitarily equivalent to the group representation
$\zeta: = \sigma_{1}\oplus\sigma_{2}$ acting on the Hilbert space
$K: = M\oplus M^{\perp}$. Since $\sigma_{1}$ is unitarily
equivalent to the right regular representation of $G$ (because of the Riesz sequence), we have
that $\sigma_{1}(G)'' \simeq \lambda(G)'$. Let $q$ be the
orthogonal projection from $K$ onto $M\oplus 0$. Then $q\in
\zeta(G)'$. Clearly, $\zeta(G)''q\simeq \sigma_{1}(G)''$. Since
$\zeta(G)''$ is a factor, we also have that $\zeta(G)''\simeq
\zeta(G)''q$, and hence $\sigma(G)'' \simeq \lambda(G)'$, i.e.,
 $\lambda(G)'\simeq P\lambda(G)'P$ since
 $\sigma(G)''=P\lambda(G)'P$.
\qed

\begin{rem}\label{remak1} Although we stated the result in Theorem
\ref{thm-ICC} for group representations, the proof works for
general projective unitary representations when the von Neumann
algebra generated by the left regular projective unitary
representation of $G$ is a factor.
\end{rem}

\noindent{\bf Proof of Theorem \ref{thm-main2}}

(i) The abelian group case is proved in Proposition
\ref{abeian-thm-main}. If $G$ is an amenable ICC group $G$, then
the statement follows immediately from Theorem \ref{thm-ICC} and
the famous result of A. Connes \cite{Conn} that when $G$ is an
amenable ICC group, then $\lambda(G)'$ is the hyperfinite $\rm
II_1$ factor, and we have that $\lambda(G)'$ and $P\lambda(G)'P$
are isomorphic for any non-zero projection $P\in \lambda(G)'$.

(ii) Recall that the fundamental group of a type $\rm II_1$ factor
$\mathcal{M}$ is the set of numbers $t>0$ for which the
''amplification'' of $\mathcal{M}$ by $t$ is isomorphic to
$\mathcal{M}$, $F(\mathcal{M})= \{t > 0: \mathcal{M}\simeq
\mathcal{M}^{t}\}$. Let $G =\Z^{2}\rtimes SL(2, \Z)$. Then, by
\cite{Popa_PNSA, Popa_AnnMath}, the fundamental group of
$\lambda(G)'$ is $\{1\}$, which implies that von Neumann algebras
$P\lambda(G)'$ is not $*$-isomorphic to $\lambda(G)'$ for any
nontrivial projection $P\in \lambda(G)'$. Thus, by Theorem
\ref{thm-ICC}, none of the nontrivial subrepresentations
$\lambda|_{P}$ admits a dual pair.\qed
\bigskip


\begin{exam}\label{exam3} Let $G = \mathcal{F}_{\infty}$. Using
Voiculescu's free probability theory, Radulescu \cite{Rad-Inv}
proved that fundamental group $F(\mathcal{M}) = \R_{+}^{*}$ for
$\mathcal{M} = \lambda(\mathcal{F}_{\infty})'$. Therefore for
$\lambda(\mathcal{F}_{\infty})'\simeq
P\lambda(\mathcal{F}_{\infty})'P$ for any nonzero projection $P\in
\lambda(\mathcal{F}_{\infty})'$, and thus $\lambda|_{P}$ admits a
dual pair for free group $\mathcal{F}_{\infty}$.
\end{exam}

\begin{exam}\label{exam4} Let $G = \Z^{d}\times \Z^{d}$, and $\pi_{\Lambda}(m, n) =
E_{m}T_{n}$ be the Gabor representation of $G$ on $L^{2}(\R^{d})$
associated with the time-frequency lattice $\Lambda= A\Z^{d}\times
B\Z^{d}$. Since $G$ is abelian, we have that the von Neumann
algebra  $\pi_{\Lambda}(G)'$ is amenable (cf. \cite{BC-08}). Thus,
if $\pi_{\Lambda}(G)'$ is a factor, then for every $\pi_{\Lambda}$
invariant subspace $M$ of $L^{2}(\R^d)$, we have by the remark after
the proof of Theorem \ref{thm-ICC} that $\pi_{\Lambda}|_{M}$
admits a dual pair. Therefore the duality principle in Gabor
analysis holds also for subspaces at least for the factor case (e.g.,
$d=1$, $A=a$ and $B=b$ with $ab$ irrational). In the case that $A
= B = I$, then the Gabor representation $\pi_{\Lambda}$ is a
unitary representation of the abelian group $\Z^{d}\times \Z^{d}$,
and so, from Proposition \ref{abeian-thm-main}, the duality
principle holds for subspaces for this case as well, In fact in
this case a concrete representation $\sigma$ can be constructed by
using the Zak transform.
\end{exam}

\begin{acknowledgements}
The authors  thank their colleagues Ionut Chifan and Palle Jorgensen for many
helpful conversations and comments on this paper.
\end{acknowledgements}

%
%
%

\end{document}